\documentclass[10pt]{article}
\usepackage{mathrsfs}
\usepackage{amsthm}
\usepackage{amsmath}
\usepackage{graphicx}
\usepackage{color}
\usepackage{amsfonts}
\newtheorem{theorem}{Theorem}[section]

\newtheorem{lemma}[theorem]{Lemma}

\numberwithin{equation}{section}
\normalsize

\begin{document}
\title{\textbf{Critical infection rates for contact processes on open clusters of oriented percolation in $Z^d$}}

\author{Xiaofeng Xue \thanks{\textbf{E-mail}: xuexiaofeng@ucas.ac.cn \textbf{Address}: School of Mathematical Sciences, University of Chinese Academy of Sciences, Beijing 100049, China.}\\ University of Chinese Academy of Sciences}

\date{}
\maketitle

\noindent {\bf Abstract}

In this paper we are concerned with contact processes on open clusters of oriented percolation in $Z^d$, where the disease spreads along the direction of open edges. We show that the two critical infection rates in the quenched and annealed cases are equal with probability one and are asymptotically equal to $(dp)^{-1}$ as the dimension $d$ grows to infinity, where $p$ is the probability of edge `open'.


\section{Introduction}
In this paper we are concerned with contact processes on open clusters of oriented bond percolation in $Z^d$. In our model, for any $x,y\in Z^d$, there is a directed edge
from $x$ to $y$ if and only if $y-x\in \{e_i\}_{1\leq i\leq d}$, where
\[
e_i=(0,~\ldots,~0,\mathop 1\limits_{i \text{th}},0,~\ldots,~0).
\]
We denote by $E_d$ the set of directed edges on $Z^d$. $\{X_e\}_{e\in E_d}$ are independent and identically distributed random variables such that
\[
P(X_e=1)=1-P(X_e=0)=p\in (0,1).
\]
Edge $e$ is called `open' if $X_e=1$ or `closed' else if $X_e=0$. We denote by $x\rightarrow y$ when the edge from $x$ to $y$ is open. After deleting all the closed edges, we obtain an oriented subgraph $G$ of $Z^d$, which our contact process will be defined on. Please note that $G$ is a random graph depending on the values of $\{X_e\}_{e\in E_d}$.

Contact processes $\{\eta_t\}_{t\geq 0}$ on $G$ is a spin system with state space $\{0,1\}^G$, which means that at each vertex, there is a spin with value $0$ or $1$. The flip rates of $\{\eta_t\}_{t\geq 0}$ are given by
\[
c(x,\eta)=
\begin{cases}
1&\text{~if~}\eta(x)=1,\\
\lambda\sum_{y:y\rightarrow x}\eta(y)&\text{~if~} \eta(x)=0
\end{cases}
\]
for any $(x,\eta)\in G\times \{0,1\}^G$, where $\lambda>0$ is the infection rate. For more details about spin systems, please see Chapter 3 of \cite{LIG1985}.

Intuitively, the model describes the spread of an infectious disease along the direction of open edges. $1$ and $0$ represent the state `infected' and `healthy' respectively. An infected vertex waits for an exponential time with rate one to recover. For any healthy vertex $x$, if $y$ is infected and the edge from $y$ to $x$ is open, then $y$ infects $x$ at rate $\lambda$.

In real life, diseases spreading along one direction are those traveling by rivers, such as dysentery, cholera, typhoid and so on. Closed edges represent the river courses which are too dry to carry the disease.

Recently, contact processes in random environments such as percolation model is an popular topic. Here we list some results in this field which inspire us. In \cite{Chen2009} and \cite{Yao2012}, Chen and Yao prove that the complete convergence theorem holds for contact processes in two kinds of random environments on $Z^d\times Z^+$, one of which is the percolation model. In \cite{Chat2009}, Chatterjee and Durrett show that contact processes on random graphs with power law degree distribution have critical value $0$, which is not consistent with the estimation given by non-rigorous mean field approach. In \cite{Pet2011}, Peterson shows that the critical value of contact processes on complete graphs with random vertex-dependent infection rates is inversely proportional to the second moment of the weight of a vertex.

\section{Main results}
We need introduce some notations before stating the main problem we concerned with. In later sections, we denote by $P^G_\lambda$ the probability measure of the contact process with infection rate $\lambda$ on a given graph $G$, which is called the quenched measure. We denote by $E^G_\lambda$ the expectation with respect to $P^G_\lambda$. Note that $G$ depends on the values of $\{X_e\}_{e\in E^d}$, which leads to following notations. We assume that $\{X_e\}_{e\in E^d}$ are defined on the product measurable space $(\{0,1\}^{E_d},\mathcal{F}_d,\mathbb{P}_{d,p})$ (see Section 1.3 of \cite{Grim1999}), where $p$ is the probability of `open'. We denote by $\mathbb{E}_{d,p}$ the expectation with respect to $\mathbb{P}_{d,p}$. For any $\omega\in \{0,1\}^{E_d}$, we denote by $G(\omega)$ the random graph of oriented percolation depending on $\{X_e(\omega)\}_{e\in E_d}$.
We define
 \[
 P_{\lambda,d,p}(\cdot)=\mathbb{E}_{d,p}\big[P^{G(\omega)}_\lambda(\cdot)\big],
 \]
which is called the annealed measure. We denote by $E_{\lambda,d,p}$ the expectation with respect to $P_{\lambda,d,p}$.

In later sections, we write $\eta_t$ as $\eta_t^A$ when
\[
\{x:\eta_0(x)=1\}=A.
\]
If all the vertices are infected at the beginning, then we omit the superscript.

Since the contact process is attractive (see the definition of attractive in Chapter 3 of \cite{LIG1985}), it is easy to see that $P^G_\lambda(\eta_t(x)=1)$ is decreasing with $t$ for any $x\in G$ and so does $P_{\lambda,d,p}(\eta_t(\textbf{0})=1)$, where $\textbf{0}$ is the origin of $Z^d$. Furthermore, according to the basic coupling of spin systems (See Chapter 3 of \cite{LIG1985}), if $\lambda_1\geq \lambda_2$, then
\[
P_{\lambda_1}^G(\eta_t(x)=1)\geq P_{\lambda_2}^G(\eta_t(x)=1)
\]
and
\[
P_{\lambda_1,d,p}(\eta_t(\textbf{0})=1)\geq P_{\lambda_2,d,p}(\eta_t(\textbf{0})=1).
\]
As a result, the definitions of the following critical values are reasonable. For $d\geq 1$, $p\in (0,1)$ and random graph $G$ with respect to $\{X_e\}_{e\in E^d}$, we define
\begin{equation}\label{equ of def of annealed critical infection}
\lambda_c(d,p)=\sup\{\lambda:~\lim_{t\rightarrow+\infty}P_{\lambda,d,p}(\eta_t(\textbf{0})=1)=0\}
\end{equation}
and
\begin{equation}\label{equ of def of quenched critical infection}
\widehat{\lambda}_c(G)=\sup\{\lambda:~\forall x\in G, \lim_{t\rightarrow+\infty}P_\lambda^G(\eta_t(x)=1)=0\}.
\end{equation}

According to the translation invariance of our model, $P_{\lambda,d,p}(\eta_t(x)=1)$ does not reply on the choose of $x$. However, the contact process on a given $G$ is not symmetric for each vertex, which explains the difference between the two definitions.

The main problem we concerned with is the estimation of $\lambda_c(d,p)$ and $\widehat{\lambda}_c(G)$. The following theorem is our main result.

\begin{theorem}\label{theorem of main}
$(\romannumeral1)$ For any $d\geq 1$ and $p\in (0,1)$, there exist $\widehat{\lambda}_c(d,p)\geq 0$ and $A_{d,p}\in \mathcal{F}_d$ such that
\[
\mathbb{P}_{d,p}(A_{d,p})=1
\]
and
\[
\widehat{\lambda}_c(G(\omega))=\widehat{\lambda}_c(d,p)
\]
for any $\omega\in A_{d,p}$.

$(\romannumeral2)$
\begin{equation}\label{equ of annealed cri is quench cri}
\lambda_c(d,p)=\widehat{\lambda}_c(d,p).
\end{equation}

$(\romannumeral3)$ For any $p\in (0,1)$,
\begin{equation}\label{equ of 1/dp}
\lim_{d\rightarrow+\infty}dp\lambda_c(d,p)=1.
\end{equation}
\end{theorem}

Theorem \ref{theorem of main} shows that the two critical infection rates in \eqref{equ of def of annealed critical infection} and \eqref{equ of def of quenched critical infection} are equal with probability one. Furthermore, as $d$ grows to infinity, these critical infection rates are asymptotically equal to $1/(dp)$, which is inversely proportional to the expectation of open edges from a fixed vertex.

Critical infection rates for contact processes on some other graphs have similar asymptotic behaviors with that in \eqref{equ of 1/dp}. In \cite{Grif1983}, Griffeath shows that $\lambda_c\approx 1/(2d)$ for contact processes on $Z^d$. In \cite{Pem1992}, Pemantle  shows that $\lambda_c\approx 1/n$ for contact process on regular tree $T^n$. In \cite{Pet2011}, Peterson shows that $\lambda_c\approx \frac{1}{nE\rho^2}$ for contact process on complete graph $C_n$ with random vertex-dependent infection rate $\rho(\cdot)$. All these results including \eqref{equ of 1/dp} are consistent with the the non-rigorous mean field analysis of contact processes. However, Chatterjee and Durrett prove in \cite{Chat2009} that contact processes on random graphs with power low degree distribution have critical infection rate $0$, hence the mean field analysis gives an incorrect estimation of critical value when the power $\alpha\geq 3$.

The proof of \eqref{equ of 1/dp} will be divided into Section \ref{section of lower bound} and Section \ref{section of upper bound}. Now we give the proof of $(\romannumeral1)$ and $(\romannumeral2)$.

\proof[Proof of $(\romannumeral1)$]

For an edge $e\in E_d$ from $x_0$ to $y_0$ and any $x\in Z^d$, we denote by $x+e$ the edge from $x+x_0$ to $x+y_0$. For any $x\in Z^d$, we define $T_x:\{0,1\}^{E_d}\rightarrow \{0,1\}^{E_d}$ as
\[
[T_x(\omega)](e)=\omega(x+e)
\]
for any $\omega\in \{0,1\}^{E_d}$ and $e\in E_d$.

It is obviously that
\[
\widehat{\lambda}_c(G(\omega))=\widehat{\lambda}_c(G(T_x(\omega)))
\]
for any $\omega\in \{0,1\}^{E_d}$ and $x\in Z^d$.

As a result, $(\romannumeral1)$ follows the ergodicity of i.i.d. measures (see Chapter 7 of \cite{Dur2010}).

\qed

\proof[Proof of \eqref{equ of annealed cri is quench cri}]

For any $\lambda<\widehat{\lambda}_c(d,p)$,
\begin{align*}
\lim_{t\rightarrow+\infty}P_{\lambda,d,p}(\eta_t(\textbf{0})=1)
=&\lim_{t\rightarrow+\infty}\mathbb{E}_{d,p}P_\lambda^{G(\omega)}(\eta_t(\textbf{0})=1)\\
=&\mathbb{E}_{d,p}\lim_{t\rightarrow+\infty}P_\lambda^{G(\omega)}(\eta_t(\textbf{0})=1)\\
=&\mathbb{E}_{d,p}\big[1_{\{\omega\in A_{d,p}\}}\lim_{t\rightarrow+\infty}P_\lambda^{G(\omega)}(\eta_t(\textbf{0})=1)\big]\\
=&0
\end{align*}
according to $(\romannumeral1)$. Therefore,
\begin{equation}\label{equ of quench cri no big}
\widehat{\lambda}_c(d,p)\leq \lambda_c(d,p).
\end{equation}
For any $\lambda<\lambda_c(d,p)$ and each $x\in Z^d$,
\[
\lim_{t\rightarrow+\infty}P_{\lambda,d,p}(\eta_t(x)=1)=0.
\]
Therefore, for any $x\in Z^d$,
\[
\mathbb{E}_{d,p}\lim_{t\rightarrow+\infty}P_\lambda^{G(\omega)}(\eta_t(x)=1)=\lim_{t\rightarrow+\infty}P_{\lambda,d,p}(\eta_t(x)=1)=0.
\]
Hence, with probability one,
\[
\lim_{t\rightarrow+\infty}P_\lambda^G(\eta_t(x)=1)=0.
\]
Notice that there are countable vertices on $Z^d$. As a result, there exists $B_{d,p}\in \mathcal{F}_d$ such that
\[
\mathbb{P}_{d,p}(B_{d,p})=1
\]
and
\[
\forall x\in Z^d, \lim_{t\rightarrow+\infty}P_\lambda^{G(\omega)}(\eta_t(x)=1)=0
\]
for any $\omega\in B_{d,p}$.

We choose $\omega_0\in A_{d,p}\cap B_{d,p}$, then
\[
\lambda\leq \widehat{\lambda}_c(G(\omega_0))=\widehat{\lambda}_c(d,p).
\]
Therefore,
\begin{equation}\label{equ of annealed cri no big}
\lambda_c(d,p)\leq \widehat{\lambda}_c(d,p).
\end{equation}
\eqref{equ of annealed cri is quench cri} follows \eqref{equ of quench cri no big} and \eqref{equ of annealed cri no big}.

\qed

\section{Mean field estimation}\label{section of mean field}
In this section we utilize the mean field approach to give a non-rigorous explanation of why $\lambda_c\approx1/(dp)$. The rigorous proof will be given in Section \ref{section of lower bound} and Section \ref{section of upper bound}.

According to Hille-Yosida Theorem,
\begin{align*}
\frac{d}{dt}P_{\lambda,d,p}(\eta_t(x)=1)=&-P_{\lambda,d,p}(\eta_t(x)=1)\\
&+\lambda\sum_{i=1}^dP_{\lambda,d,p}\big(\eta_t(x)=0,\eta_t(x-e_i)=1,x-e_i\rightarrow x\big).
\end{align*}

In the mean field approach, we assume that $\eta_t(x)$, $\eta_t(x-e_i)$ and $1_{\{x-e_i\rightarrow x\}}$ are independent (which is wrong).

Then,
\begin{align*}
&P_{\lambda,d,p}\big(\eta_t(x)=0,\eta_t(x-e_i)=1,x-e_i\rightarrow x\big)\\
=&pP_{\lambda,d,p}(\eta_t(x)=1)\big(1-P_{\lambda,d,p}(\eta_t(x)=1)\big)
\end{align*}
under the mean field assumption.

So $P_{\lambda,d,p}(\eta_t(x)=1)$ is described by the following ODE
\begin{equation}\label{equ of mean field ODE}
\begin{cases}
&\frac{d}{dt}f_t=-f_t+\lambda dpf_t(1-f_t),\\
&f_0=1.
\end{cases}
\end{equation}
By direct calculation,
\[
\lim_{t\rightarrow+\infty}f_t=0
\]
when $\lambda<1/(dp)$ and
\[
\lim_{t\rightarrow+\infty} f_t=\frac{\lambda dp-1}{\lambda dp}>0
\]
when $\lambda>1/(dp)$.

As a result, the estimation of $\lambda_c(d,p)$ given by the mean field approach is $1/(dp)$, which is actually a lower bound shown in the next section.

\section{Lower bound of $\lambda_c(d,p)$}\label{section of lower bound}
In this section we give a lower bound of $\lambda_c(d,p)$. We utilize the binary contact path process as an auxiliary process, which is introduced by Griffeath in \cite{Grif1983}.

For any $\omega\in \{0,1\}^{E_d}$, the binary contact path process $\{\zeta_t\}_{t\geq0}$ on $G(\omega)$ is with state space $\{0,1,2,\ldots\}^{G(\omega)}$, which means that each vertex takes a value from nonnegative integers. $\{\zeta_t\}_{t\geq 0}$ evolves as follows. For each $x\in Z^d$, $\zeta_t(x)$ flips to $0$ at rate one. For each $y$ such that $y\rightarrow x$, $\zeta_t(x)$ flips to $\zeta_t(x)+\zeta_t(y)$ at rate $\lambda$.

In other words, the generator $\Omega$ of $\{\zeta_t\}_{t\geq 0}$ is given by
\[
\Omega f(\zeta)=\sum_{x\in Z^d}\big[f(\zeta^{x,0})-f(\zeta)\big]+\lambda\sum_{x\in Z^d}\sum_{y:y\rightarrow x}\big[f(\zeta^{x,\zeta(x)+\zeta(y)})-f(\zeta)\big]
\]
for any $\zeta\in \{0,1,2,\ldots\}^G$, where
\[
\zeta^{x,m}(y)=
\begin{cases}
\zeta(y)& \text{~if~}y\neq x,\\
m &\text{~if~}y=x
\end{cases}
\]
for any $(x,m)\in Z^d\times \{0,1,2,\ldots\}$.

Intuitively, the binary contact path process $\{\zeta_t\}_{t\geq 0}$ counts the seriousness of the disease. An infected vertex $x$ is able to be further infected by $y$ if there is an open edge from $y$ to $x$. When $y$ infects $x$, we add the seriousness of the disease of $x$ by the seriousness of $y$.

We assume that $\zeta_0(x)=1$ for any $x\in Z^d$. Then, it is easy to see that the contact process $\{\eta_t\}_{t\geq 0}$ with
\[
\{x:\eta_0(x)=1\}=Z^d
\]
can be coupled with $\{\zeta_t\}$ as follows. For any $x\in Z^d$ and $t\geq 0$,
\begin{equation}\label{equ of couple of zeta eta}
\eta_t(x)=
\begin{cases}
1 & \text{~if~}\zeta_t(x)\geq 1,\\
0 & \text{~if~}\zeta_t(x)=0.
\end{cases}
\end{equation}

By \eqref{equ of couple of zeta eta},
\[
P_\lambda^{G(\omega)}(\eta_t(\textbf{0})=1)=P_\lambda^{G(\omega)}(\zeta_t(\textbf{0})\geq 1)\leq E_\lambda^{G(\omega)}\zeta_t(\textbf{0})
\]
and hence
\begin{equation}\label{equ of Peta small than Ezeta}
P_{\lambda,d,p}(\eta_t(\textbf{0})=1)\leq E_{\lambda,d,p}\zeta_t(\textbf{0}).
\end{equation}
The following theorem gives a lower bound of $\lambda_c(d,p)$.
\begin{theorem}\label{theorem of lower bound}
For any $d\geq 1$ and $p\in (0,1)$,
\begin{equation}\label{equ of lower bound}
\lambda_c(d,p)\geq 1/(dp).
\end{equation}
\end{theorem}
\proof

According to Theorem 1.27 of Chapter 9 of \cite{LIG1985}, for any $x\in Z^d$ and any $\omega\in \{0,1\}^{E_d}$,
\[
\frac{d}{dt}E_\lambda^{G(\omega)}\zeta_t(x)=-E_\lambda^{G(\omega)}\zeta_t(x)+\lambda\sum_{y:y\rightarrow x}E_\lambda^{G(\omega)}\zeta_t(y).
\]
Therefore,
\[
E_\lambda^{G(\omega)}\zeta_t=e^{t(B_\omega-I_d)}\zeta_0,
\]
where $B_\omega$ is a $Z^d\times Z^d$ matrix such that
\[
B_\omega(x,y)=
\begin{cases}
\lambda &\text{~if~}y\rightarrow x,\\
0 & \text{~else}
\end{cases}
\]
for any $(x,y)\in Z^d\times Z^d$ and $I_d$ is the $Z^d\times Z^d$ identity matrix.

As a result,
\begin{align*}
E_\lambda^{G(\omega)}\zeta_t(\textbf{0})&=e^{-t}\sum_{y\in Z^d}e^{tB_\omega}(0,y)\\
&=e^{-t}\sum_{n=0}^{+\infty}\frac{t^n}{n!}\big[\sum_{y\in Z^d}B_\omega^n(0,y)\big].
\end{align*}
For any $\omega\in \{0,1\}^{E_d}$ and $n\geq 0$, we denote by $l_n(\omega)$ the total number of open paths which are with length $n$ and end at $\textbf{0}$.

By the definition of $B_\omega$, it is easy to see that
\[
\sum_{y\in Z^d}B_\omega^n(0,y)=\lambda^nl_n(\omega).
\]
Therefore,
\[
E_{\lambda,d,p}\zeta_t(\textbf{0})=e^{-t}\sum_{n=0}^{+\infty}\frac{t^n\lambda^n}{n!}\mathbb{E}_{d,p}l_n.
\]
For the oriented percolation in $Z^d$, there are $d^n$ paths to $\textbf{0}$ with length $n$. Each path is open with probability $p^n$.
Hence,
\[
\mathbb{E}_{d,p}l_n=d^np^n
\]
and
\[
E_{\lambda,d,p}\zeta_t(\textbf{0})=e^{-t}\sum_{n=0}^{+\infty}\frac{t^n\lambda^nd^np^n}{n!}=\exp\{(\lambda dp-1)t\}.
\]
Therefore, when $\lambda<1/(dp)$,
\begin{equation}\label{equ of exponential decay of zeta}
\lim_{t\rightarrow+\infty}E_{\lambda,d,p}\zeta_t(\textbf{0})=\lim_{t\rightarrow+\infty}\exp\{(\lambda dp-1)t\}=0.
\end{equation}
\eqref{equ of lower bound} follows \eqref{equ of Peta small than Ezeta} and \eqref{equ of exponential decay of zeta}.

\qed

As a direct corollary of Theorem \ref{theorem of lower bound},
\begin{equation}\label{equ of liminf cri}
\liminf_{d\rightarrow+\infty}dp\lambda_c(d,p)\geq 1.
\end{equation}

\section{Upper bound of $\lambda_c(d,p)$}\label{section of upper bound}
In this section we give an upper bound of $\lambda_c(d,p)$. We are inspired a lot by the approaches in \cite{Cox1983} and \cite{Kesten1990}, which are introduced by Kesten.

According to the graphic representation of contact processes (see Section 3.6 of \cite{LIG1985}), $\{\eta_t\}_{t\geq 0}$ has a dual process $\{\widehat{\eta}_t\}_{t\geq 0}$, where the disease spreads along the opposite direction of the edge. In details, $\{\widehat{\eta}_t\}_{t\geq 0}$ on graph $G$ is a spin system with flip rates function given by
\[
\widehat{c}(x,\eta)=
\begin{cases}
1 &\text{~if~}\eta(x)=1,\\
\lambda\sum_{y:x\rightarrow y}\eta(y) &\text{~if~}\eta(x)=0.\\
\end{cases}
\]
We write $\widehat{\eta}_t$ as $\widehat{\eta}_t^A$ when $\{x:\widehat{\eta}_0(x)=1\}=A$. The graphic representation shows that
\[
P_\lambda^G(\eta_t(\textbf{0})=1)=P_\lambda^G(\exists ~x\in Z^d, \widehat{\eta}_t^\textbf{0}(x)=1)
\]
and hence
\[
P_{\lambda,d,p}(\eta_t(\textbf{0})=1)=P_{\lambda,d,p}(\exists ~x\in Z^d, \widehat{\eta}_t^\textbf{0}(x)=1).
\]
It is easy to see that $\{\exists ~x\in Z^d, \widehat{\eta}_t^\textbf{0}(x)=1\}$ and $\{\exists ~x\in Z^d, \eta_t^\textbf{0}(x)=1\}$ has the same distribution under the annealed measure $P_{\lambda,d,p}$. As a result,
\begin{equation}\label{equ of self-dual}
P_{\lambda,d,p}(\eta_t(\textbf{0})=1)=P_{\lambda,d,p}(\exists ~x\in Z^d, \eta_t^\textbf{0}(x)=1).
\end{equation}

\eqref{equ of self-dual} gives the self-duality of $\eta_t$, but please note that in our model this self-duality only holds in the annealed case, not in the quenched case.

By \eqref{equ of self-dual},
\begin{equation}\label{equ of survives}
\lim_{t\rightarrow+\infty}P_{\lambda,d,p}(\eta_t(\textbf{0})=1)=P_{\lambda,d,p}(\eta_t^{\textbf{0}} \text{~survives}).
\end{equation}

We control the evolution of $\eta_t^\textbf{0}$ from below by a SIR model. Assume that $\{Y_x\}_{x\in Z^d}$ and $\{U_{x,i}\}_{x\in Z^d, 1\leq i\leq d}$ are independent random variables and are independent with $\{X_e\}_{e\in E_d}$. For any $(x,i)\in Z^d\times \{1,2,\ldots,d\}$, $Y_x$ follows exponential distribution with rate one and $U_{x,i}$ follows exponential distribution with rate $\lambda$. For any $(x,i)\in Z^d\times \{1,2,\ldots,d\}$, if the edge from $x$ to $y=x+e_i$ is open and $U_{x,i}\leq Y_x$, then we say that $x$ infects $y$, which is denoted by $x\Rightarrow y$.

For any $x\in Z^d$ such that $\sum_{i=1}^dx(i)=n$, if there exists $\{y_i:y_i\in Z^d\}_{1\leq i\leq n-1}$ such that $\textbf{0}\Rightarrow y_1$, $y_i\Rightarrow y_{i+1}$ for $i\leq n-2$ and $y_{n-1}\Rightarrow x$, then we say that there is an infection path with length $n$ from $\textbf{0}$ to $x$.

We denote by $C_n$ the set of infection pathes from $\textbf{0}$ with length $n$. It is easy to see that
\[
\{\eta_t^{\textbf{0}}\text{~survives}\}\supseteq \{\forall ~n\geq1, C_n\neq \emptyset\}
\]
in the sense of coupling.

As a result,
\begin{equation}\label{equ of SIR small than SI}
P_{\lambda,d,p}(\eta_t^{\textbf{0}} \text{~survives})\geq \lim_{n\rightarrow+\infty}P_{\lambda,d,p}(C_n\neq \emptyset).
\end{equation}

To give lower bound of $P_{\lambda,d,p}(C_n\neq \emptyset)$, we utilize the processes of simple random walk on oriented lattices. Assume that $\{S_m^d\}_{m=0,1,2,\ldots}$
is a random walk on $Z^d$ such that $S_0^d=\textbf{0}$ and
\[
P(S_{m+1}^d-S_m^d=e_i)=1/d
\]
for $1\leq i\leq d$. Let $\{\widehat{S}_m^d\}_{m=0,1,2,\ldots}$ be an independent copy of $\{S_m^d\}_{m=0,1,2,\ldots}$.

We define
\[
k_d=\sum_{i=0}^{+\infty}1_{\{S_i^d=\widehat{S}_i^d,S_{i+1}^d\neq \widehat{S}_{i+1}^d\}}
\]
and
\[
r_d=\sum_{i=0}^{+\infty}1_{\{S_i^d=\widehat{S}_i^d,S_{i+1}^d=\widehat{S}_{i+1}^d\}}.
\]
The following lemma is crucial for us to estimate $\lambda_c(d,p)$.
\begin{lemma}\label{lemma of upper bound}
If $\lambda$ satisfies that
\[
E\big[2^{k_d}(\frac{\lambda+1}{\lambda p})^{r_d}\big]<+\infty,
\]
then
\[
\lambda\geq \lambda_c(d,p).
\]
\end{lemma}

\proof

For any $n\geq 1$, we define
\[
T_n=\Big\{\{x_k\}_{k=0}^n:x_0=\textbf{0},x_{j+1}-x_j\in \{e_i:1\leq i\leq d\}\text{~for~}0\leq j\leq n-1\Big\}
\]
as the set of pathes from $\textbf{0}$ with length $n$ (no matter whether each edge is open or closed).

Let
\[
M_n=\Big\{\{x_k\}_{k=0}^n\in T^n:x_i\Rightarrow x_{i+1}\text{~for~}0\leq i\leq n-1\Big\}\subseteq T_n.
\]
Then,
\begin{equation}\label{equ of Cn holder}
P_{\lambda,d,p}(C_n\neq\emptyset)=P(|M_n|>0)\geq \frac{(E|M_n|)^2}{E|M_n|^2}.
\end{equation}
For $x,y_1,y_2$ such that $y_1-x,y_2-x\in \{e_j:1\leq j\leq d\}$,
\begin{equation}\label{equ of infect one}
P(x\Rightarrow y_1)=pP(U_{x,1}\leq Y_x)=\frac{\lambda p}{1+\lambda}
\end{equation}
and
\begin{equation}\label{equ of infect two}
P(x\Rightarrow y_1, ~x\Rightarrow y_2)=p^2P(U_{x,1},U_{x,2}\leq Y_x)=\frac{2\lambda^2p^2}{(2\lambda+1)(\lambda+1)}.
\end{equation}

Let
\[
A_n=\{0\leq i\leq n-1:S^d_i=\widehat{S}^d_i,S^d_{i+1}=\widehat{S}^d_{i+1}\}
\]
and
\[
B_n=\{0\leq i\leq n-1:S^d_i=\widehat{S}^d_i,S^d_{i+1}\neq\widehat{S}^d_{i+1}\}.
\]

Then by \eqref{equ of Cn holder}, \eqref{equ of infect one} and \eqref{equ of infect two},
\begin{align}\label{equ of lower bound of PCn}
P_{\lambda,d,p}(C_n\neq\emptyset)&\geq\frac{d^{2n}(\frac{\lambda p}{\lambda+1})^{2n}}{\sum\limits_{\{y_i\}_{i=0}^n\in T_n}\sum\limits_{\{z_i\}_{i=0}^n\in T_n}
P\big(\{y_i\}_{i=0}^n\in M_n,\{z_i\}_{i=0}^n\in M_n\big)}\notag\\
&=\frac{(\frac{\lambda p}{1+\lambda})^{2n}}{P\big(\{S_i^d\}_{i=0}^n\in M_n, \{\widehat{S}_i^d\}_{i=0}^n\in M_n\big)}\notag\\
&=\frac{(\frac{\lambda p}{1+\lambda})^{2n}}{E\Big[(\frac{\lambda p}{1+\lambda})^{2n-|A_n|-2|B_n|}(\frac{2\lambda^2 p^2}{(1+2\lambda)(1+\lambda)})^{|B_n|}\Big]}\notag\\
&=\frac{1}{E\Big[(\frac{1+\lambda}{\lambda p})^{|A_n|}(\frac{2+2\lambda}{1+2\lambda})^{|B_n|}\Big]}\notag\\
&\geq \frac{1}{E\Big[2^{|B_n|}(\frac{1+\lambda}{\lambda p})^{|A_n|}\Big]}.
\end{align}

Notice that $\lim_{n\rightarrow+\infty}|A_n|=r_d$ and $\lim_{n\rightarrow+\infty}|B_n|=k_d$. Therefore, by \eqref{equ of survives}, \eqref{equ of SIR small than SI} and \eqref{equ of lower bound of PCn},
\[
\lim_{t\rightarrow+\infty}P_{\lambda,d,p}(\eta_t(\textbf{0})=1)\geq \frac{1}{E\big[2^{k_d}(\frac{\lambda+1}{\lambda p})^{r_d}\big]}>0
\]
when
\[
E\big[2^{k_d}(\frac{\lambda+1}{\lambda p})^{r_d}\big]<+\infty.
\]

\qed

Thanks to Lemma \ref{lemma of upper bound}, we can give an upper bound of $\lambda_c(d,p)$ to finish the proof of Theorem \ref{theorem of main}. 

\proof[Proof of $~\limsup\limits_{d\rightarrow+\infty}dp\lambda_c(d,p)\leq 1$.]

\quad

We define a sequence of increasing stopping times $\{\tau_n\}_{n=1}^{+\infty}$ about $\{S_n^d, \widehat{S}_n^d\}_{n=0,1,2,\ldots,}$ as follows.
\[
\tau_1=\inf\{k\geq 0:S^d_k=\widehat{S}^d_k,S^d_{k+1}=\widehat{S}^d_{k+1}\}.
\]
If $\tau_1=+\infty$, then $\tau_j=+\infty$ for $j\geq 2$. If $\tau_1<\infty$,
then
\[
\tau_2=\inf\{k>\tau_1:S^d_k=\widehat{S}^d_k,S^d_{k+1}=\widehat{S}^d_{k+1}\}.
\]
By inducing, if $\tau_l=+\infty$, then $\tau_j=+\infty$ for $j\geq l+1$.
If $\tau_l<\infty$, then
\[
\tau_{l+1}=\inf\{k>\tau_l:S^d_k=\widehat{S}^d_k,S^d_{k+1}=\widehat{S}^d_{k+1}\}.
\]
We set $\tau_0=-1$ for later use.
For $k\geq1$, if $\tau_k<+\infty$, then we define
\[
\sigma_k=
\begin{cases}
0& \text{~if~}\tau_k=\tau_{k-1}+1,\\
\sum_{j=\tau_{k-1}+1}^{\tau_k-1}1_{\{S_j^d=\widehat{S}_j^d\}}&\text{~if~}\tau_k-\tau_{k-1}\geq 2.\\
\end{cases}
\]
If $\tau_k<+\infty$ and $\tau_{k+1}=+\infty$, then we define
\[
\rho_k=\sum_{j=\tau_k+1}^{+\infty}1_{\{S_j^d=\widehat{S}_j^d\}}.
\]
We define
\[
\theta=\inf\{j\geq 1:S_j^d=\widehat{S}_j^d\}.
\]
We write $\tau_k$, $\sigma_k$, $\rho_k$ and $\theta$ as $\tau_k(d)$, $\sigma_k(d)$, $\rho_k(d)$ and $\theta(d)$ when the dimension $d$ need to be distinguished.

According to Markov property,
\begin{equation}\label{equ of Esigma}
P(\sigma_k=m,\tau_k<+\infty\big|\tau_{k-1}<+\infty)=P(2\leq \theta<+\infty)^m\frac{1}{d}
\end{equation}
and
\begin{equation}\label{equ of Erho}
P(\rho_k=l,\tau_{k+1}=+\infty\big|\tau_k<+\infty)=P(2\leq \theta<+\infty)^{l-1}P(\theta=+\infty).
\end{equation}
If $r_d=k$, then $\tau_k(d)<+\infty$, $\tau_{k+1}(d)=+\infty$ and
\[
k_d=\sum_{l=1}^k\sigma_l(d)+\rho_k(d).
\]
Therefore, according to Markov property,
\begin{align}\label{equ of 4.10}
E\big[2^{k_d}(\frac{\lambda+1}{\lambda p})^{r_d}\big]
=&\sum_{k=0}^{+\infty}(\frac{\lambda+1}{\lambda p})^kE\Big[2^{^{\sum_{l=1}^k\sigma_l(d)+\rho_k(d)}}1_{\{\tau_k(d)<+\infty,\tau_{k+1}(d)=+\infty\}}\Big]\notag\\
=&\sum_{k=0}^{+\infty}(\frac{1+\lambda}{\lambda p})^k\Big(E\big[2^{\sigma_1(d)}1_{\{\tau_1(d)<+\infty\}}\big]
\Big)^kE\big[2^{\rho_0(d)}1_{\{\tau_1(d)=+\infty\}}\big].
\end{align}
It is proven by Cox and Durrett in \cite{Cox1983} that there exists $C>0$ such that
\begin{equation}\label{equ of 4.11}
P(2\leq \theta(d)<+\infty)\leq \frac{C}{d^2}
\end{equation}
for any $d\geq 1$.

By \eqref{equ of Esigma}, \eqref{equ of Erho} and \eqref{equ of 4.11}, for $d>2\sqrt{C}$,
\begin{align*}
E\big[2^{\sigma_1(d)}1_{\{\tau_1(d)<+\infty\}}\big]
=&\sum_{l=0}^{+\infty}2^lP(2\leq \theta(d))<+\infty)^l\frac{1}{d}\\
\leq&\frac{1}{d}\sum_{l=0}^{+\infty}(\frac{2C}{d^2})^l=\frac{d}{d^2-2C}
\end{align*}
and
\begin{align*}
E\big[2^{\rho_0(d)}1_{\{\tau_1(d)=+\infty\}}\big]
\leq&\sum_{l=1}^{+\infty}2^lP(2\leq \theta(d)<+\infty)^{l-1}\\
\leq&2\sum_{l=0}^{+\infty}(\frac{2C}{d^2})^l=\frac{2d^2}{d^2-2C}<4.
\end{align*}
Therefore, by \eqref{equ of 4.10},
\begin{equation}\label{equ of 4.12}
E\big[2^{k_d}(\frac{\lambda+1}{\lambda p})^{r_d}\big]\leq 4\sum_{k=0}^{+\infty}\Big[\frac{d(\lambda+1)}{(d^2-2c)\lambda p}\Big]^k
\end{equation}
for sufficiently large $d$.

By Lemma \ref{lemma of upper bound} and \eqref{equ of 4.12}, $\lambda\geq\lambda_c(d,p)$ when
\[
\frac{d(1+\lambda)}{\lambda p(d^2-2C)}<1.
\]
Therefore,
\[
\lambda_c(d,p)\leq \frac{1}{dp-\frac{2pC}{d}-1}
\]
for sufficiently large $d$ and
\[
\limsup_{d\rightarrow+\infty}dp\lambda_c(d,p)\leq 1.
\]

\qed

Since we have shown that $\liminf_{d\rightarrow+\infty}dp\lambda_c(d,p)\geq 1$ in Section \ref{section of lower bound}, the whole proof of \eqref{equ of 1/dp} and Theorem \ref{theorem of main} is completed.

{}
\end{document}